# CORRECTIONS AND ACKNOWLEDGMENT FOR "LOCAL LIMIT THEORY AND LARGE DEVIATIONS FOR SUPERCRITICAL BRANCHING PROCESSES"

By P. E. Ney and Anand N. Vidyashankar

*University of Wisconsin and Cornell University*

Theorem 1 in [2] is incorrect in the case $\alpha \geq 1$. Our error stems from the fact that the lower bound $C_1$ was determined by an integral expression which we treated as positive, whereas in fact it was zero. This led to an incorrect normalization $A_n$ when $\alpha \geq 1$. This error was communicated to us by K. Fleischmann and V. Wachtel, and the correction, that $A_n = p_1^n v_n^{(\alpha-1)}$ for all $0 < \alpha < \infty$, appears in their paper [1]. We thank them for this communication. The same error carried into Theorem 2, where the inequality (8) holds for all $0 < \alpha < \infty$.

Department of Mathematics
University of Wisconsin
Madison, Wisconsin 53706-1388
USA

Department of Statistical Science
Cornell University
Ithaca, New York 14853-4201
USA
E-mail: anv4@cornell.edu